\documentclass[11pt]{article}
\usepackage[centertags]{amsmath}
\usepackage{amsfonts}
\usepackage{amssymb}
\usepackage{amsthm}
\usepackage{newlfont}
\usepackage{graphicx}
\theoremstyle{definition}

\theoremstyle{remark}

\theoremstyle{theorem}

\begin{document}
\thispagestyle{empty}

%
%

\title{ Solving generalized Abel's integral equations of the first and second kinds via Taylor-collocation method }
\author{Eisa Zarei $^{a,}$ \footnote{Corresponding author, Tel.: +98 918 370 9587, E-mail address: ei.zarei1353@gmail.com.} ~~~~~~ \ and~~~~~~~ Samad Noeiaghdam $^{b,}$ \footnote{ E-mail addresses: s.noeiaghdam.sci@iauctb.ac.ir,
samadnoeiaghdam@gmail.com.}\
 }
  \date{}
 \maketitle

\begin{center}
\scriptsize{$^{a}$ Department of Mathematics, Hamedan
Branch, Islamic Azad University,  Hamedan, Iran. \\
$^{b}$ Department of Mathematics, Central Tehran Branch, Islamic
Azad University, Tehran, Iran.}
\end{center}
\begin{abstract}
The aim of this paper is to present an efficient numerical procedure
to approximate the generalized Abel's integral equations of the
first and second kinds. For this reason, the Taylor polynomials and
the collocation method are applied. Also, the error analysis of
presented method is illustrated. Several examples are approximated
and the numerical results show the accuracy and efficiency of this
method.

 \vspace{.5cm}{\it keywords:} Generalized Abel's integral equation, Collocation
 method, Taylor polynomials.\\


\end{abstract}
\section{Introduction}

The great mathematician Niels Abel, gave the initiative of integral
equations in 1823 in his study of mathematical physics
\cite{b5,a1,a2,a4}. Zeilon in 1924, studied the generalized Abel's
integral equation on a finite segment \cite{a5}.

The Abel's integral equations are the singular form of Volterra
integral equations. The singular integral equations are the
important and applicable kinds of integral equations which solved by
many authors \cite{myy,me1,me2,d7,b10,a4}. These integral equations
have many applications in various areas such as simultaneous dual
relations \cite{b1}, stellar winds \cite{b2}, water wave \cite{b3},
spectroscopic data \cite{b4} and the others \cite{b7,b5,b6}.

In this paper, we investigate the following generalized Abel's
integral equations of the first kind
\begin{equation}\label{0-1}
\int_a^x \frac{\Phi(t)}{(\phi(x)-\phi(t))^\alpha}dt =
g(x),~~~~~~~~0<\alpha<1,
\end{equation}
and the second kind
\begin{equation}\label{0-2}
\Phi(x)= g(x)+ \int_a^x
\frac{\Phi(t)}{(\phi(x)-\phi(t))^\alpha}dt,~~~~~~~~0<\alpha<1,
\end{equation}

where $a$ is a given real value, $g(x)$ and $\phi(x)$ are known
functions and $\Phi(x)$ is an unknown function that $\phi(t)$ is
strictly monotonically increasing and differentiable function in
some interval $a < t < b$, and $\phi'(t) \neq 0$ for every $t$ in
the interval.

Abel's integral equations (\ref{0-1}) and (\ref{0-2}) have
applications in experimental physics such as plasma diagnostics,
physical electronics, nuclear physics, optics and astrophysics
\cite{b2,b6}. Also, Baker \cite{5}, Wazwaz \cite{a1,a2} and Delves
\cite{8} studied the numerical treatment of singular integral
equations of the first and second kinds. In recent years, many
authors solved the Abel's integral equations of the first and second
kinds by using different applicable methods
\cite{d2,d1,d6,d4,d7,d3,d5}. The properties of Abel's integral
equations can be found in \cite{b11,b9,b10,a2}.

One of powerful and efficient methods to estimate the integral and
differential equations is the collocation method
\cite{col2,col1,col4}. This method is one of the expansion methods
which is used with different basis functions. In this paper, the
collocation method and the Taylor polynomials are applied to
estimate the generalized Abel's integral equations of the first and
second kinds.

Kanwal and Liu in \cite{e5} presented the Taylor expansion approach
for solving integral equations. It was applied for solving the
Volterra- Fredholm integral equations \cite{e1,e4,e7}, system of
integral equations \cite{e3}, Volterra integral equations \cite{e6},
integro-differential equations \cite{e9,e8} and the others
\cite{e2,e5}.

This paper is organized as follows: At first, in Section 2, the
Taylor polynomials and the collocation method are combined to
estimate the generalized form of Abel's integral equations of the
first and second kinds. In Section 3, the error analysis of
presented method is illustrated and in sequel, in Section 4, several
applicable examples are approximated. Also, in this section the
numerical results and the absolute errors for different values of
$x$ and $n$ are tabulated. Finally, Section 5 is conclusion.

\section{Main Idea}
In order to estimate the Abel's integral equations of the first and
second kinds, the Taylor polynomials of degree $n$ at $x=z$ is
introduced as follows
\begin{equation}\label{2-1}
\displaystyle \Phi_n(x) = \sum_{j=0}^n \frac{1}{j!} \Phi^{(j)}(z)
(x-z)^j,~~~~~~~~a\leq x,z \leq b,
\end{equation}
where the coefficients $\Phi^{(j)}(z), j=0,1,..., n$ should be
determined.

To approximate the Abel's integral equation of the first kind
(\ref{0-1}), we rewrite it in the form
\begin{equation}\label{2-2}
\displaystyle \sum_{j=0}^n \frac{1}{j!} \int_a^x \frac{
\Phi^{(j)}(z) (t-z)^j}{\big[\phi(x)-\phi(t)\big]^\alpha}dt =
g(x),~~~~~~~~0<\alpha<1~~and~~a\leq x,z \leq b.
\end{equation}
In order to estimate the integral equation (\ref{2-2}), the
collocation points
\begin{equation}\label{2-3}
x_i = a+\left(\frac{b-a}{n}\right) i,~~~~~i=0,1,2,...,n,
\end{equation}
are substitute into Eq. (\ref{2-2}) as follows
\begin{equation}\label{2-4}
\displaystyle \sum_{j=0}^n \frac{1}{j!} \int_a^{x_i} \frac{
(t-z)^j}{\big[\phi(x_i)-\phi(t)\big]^\alpha}dt~\Phi^{(j)}(z) =
g(x_i),~~~~~~~~0<\alpha<1.
\end{equation}

Now, we can write Eq. (\ref{2-4}) in the form
\begin{equation}\label{2-5}
AX=G,
\end{equation}
where
$$
A= \left[
\begin{array}{cccc}
   \frac{1}{0!} \int_a^{x_0} \frac{
(t-z)^0}{\big[\phi(x_0)-\phi(t)\big]^\alpha}dt & \frac{1}{1!}
\int_a^{x_0} \frac{ (t-z)^1}{\big[\phi(x_0)-\phi(t)\big]^\alpha}dt
&\cdots & \frac{1}{n!} \int_a^{x_0} \frac{
(t-z)^n}{\big[\phi(x_0)-\phi(t)\big]^\alpha}dt \\
   \\
 \frac{1}{0!} \int_a^{x_1} \frac{
(t-z)^0}{\big[\phi(x_1)-\phi(t)\big]^\alpha}dt & \frac{1}{1!}
\int_a^{x_1} \frac{ (t-z)^1}{\big[\phi(x_1)-\phi(t)\big]^\alpha}dt
&\cdots & \frac{1}{n!} \int_a^{x_1} \frac{
(t-z)^n}{\big[\phi(x_1)-\phi(t)\big]^\alpha}dt \\
\\
   \vdots & \vdots & \cdots & \vdots \\
   \\
   \frac{1}{0!} \int_a^{x_n} \frac{
(t-z)^0}{\big[\phi(x_n)-\phi(t)\big]^\alpha}dt & \frac{1}{1!}
\int_a^{x_n} \frac{ (t-z)^1}{\big[\phi(x_n)-\phi(t)\big]^\alpha}dt
&\cdots & \frac{1}{n!} \int_a^{x_n} \frac{
(t-z)^n}{\big[\phi(x_n)-\phi(t)\big]^\alpha}dt \\
  \end{array}
  \right]_{(n+1)(n+1)},
$$
$$
\begin{array}{cl}
  X = \left[ \begin{array}{c}
       \Phi^{(0)}(z) \\
       \\
       \Phi^{(1)}(z) \\
       \\
       \vdots \\
       \\
       \Phi^{(n)}(z)
     \end{array}\right]_{(n+1)\times 1},
     &
      G = \left[ \begin{array}{c}
       g(x_0) \\
       \\
       g(x_1) \\
       \\
       \vdots \\
       \\
       g(x_n)
     \end{array}\right]_{(n+1)\times 1}.
\end{array}
$$
By solving the system of equations (\ref{2-5}) and substitute
$\Phi^{(j)}(z), j=0,1,..., n$  in Eq. (\ref{2-1}) the approximate
solution of Abel's integral equation of the first kind (\ref{0-1})
can be obtained.

In order to estimate the second kind form of Abel's integral
equation ({\ref{0-2}}), by substitute Eq. (\ref{2-1}) into Eq.
(\ref{0-2}) we have
\begin{equation}\label{2-6}
\displaystyle \sum_{j=0}^n \frac{1}{j!} \Phi^{(j)}(z) (x-z)^j -
\int_a^x \frac{\sum_{j=0}^n \frac{1}{j!} \Phi^{(j)}(z)
(t-z)^j}{(\phi(x)-\phi(t))^\alpha}dt= g(x),~~~~~~~~0<\alpha<1,
\end{equation}
and by substitute collocation points (\ref{2-3}) into Eq.(\ref{2-6})
the following equation is obtained
\begin{equation}\label{2-7}
\sum_{j=0}^n \frac{1}{j!} \bigg[  (x_i-z)^j - \int_a^{x_i} \frac{
(t-z)^j}{(\phi(x_i)-\phi(t))^\alpha}dt \bigg] \Phi^{(j)}(c) =
g(x_i),~~~~~~~~0<\alpha<1.
\end{equation}

Finally, we rewrite Eq. (\ref{2-7}) in the matrix form $(B-A) X = G$
where matrices $A, X, G$ were presented and
$$
B= \left[
\begin{array}{cccc}
\frac{1}{0!} (x_0-z)^0 & \frac{1}{1!} (x_0-z)^1 &\cdots &
\frac{1}{n!} (x_0-z)^n \\
   \\
\frac{1}{0!}(x_1-z)^0& \frac{1}{1!} (x_1-z)^1 &\cdots & \frac{1}{n!}
(x_1-z)^n \\
\\
   \vdots & \vdots & \cdots & \vdots \\
   \\
 \frac{1}{0!} (x_n-z)^0 & \frac{1}{1!} (x_n-z)^1 &\cdots & \frac{1}{n!}
(x_n-z)^n\\
  \end{array}
  \right]_{(n+1)(n+1)}.
$$

By solving the obtained system and determine the unknowns $X$ and
substitute into Eq. (\ref{2-1}), the approximate solution of Eq.
(\ref{0-2}) can be obtained.

\section{Error analysis }
In this section, the error analysis theorem of this method is
presented.

\emph{\textbf{Theorem 1:}} Let $\Phi(x)$ is an exact solution of
Eqs. (\ref{0-1}) and (\ref{0-2}), $\Phi_n(x)$ is an approximate
solution which is obtained from $n$-th order Taylor-collocation
method and $ S_n(x) = \sum_{i=0}^n \frac{\widetilde{\Phi}^{(i)}(z)
(x-z)^i}{i !}$ is the Taylor polynomial of degree $n$ at $x=z$ then
$$
\| \Phi (x) - \Phi_n(x) \|_{\infty} \leq \frac{M}{(i+1)!}
\widetilde{\Phi}^{(i+1)}(\xi) + C \max |e_i(z)| ,~~~\xi \in [a,b],
$$
where $M = \max | (x-z)^{i+1}|$,$C = \| \ell \|_{\infty} $, $e_i(z)
= \widetilde{\Phi}^{(i)}(z) - \Phi^{(i)}(z)$.\\

\textbf{Proof:} Assume that $R_n(x)$ is the reminder term of $n$-th
order Taylor polynomial $S_n$ at $x=z$ which is given by
$$R_n(x) = \Phi (x) -S_n(x) = \frac{\widetilde{\Phi}^{(i+1)}(\xi)}{(i+1)!}
(x-z)^{i+1},\xi \in [a,b].$$
Thus
\begin{equation}\label{t0}
| R_n(x)|  = | \Phi (x) -S_n(x) | \leq
\frac{\widetilde{\Phi}^{(i+1)}(\xi)}{(i+1)!}. \max | (x-z)^{i+1}| =
\frac{M}{(i+1)!} \widetilde{\Phi}^{(i+1)}(\xi) ,~~~\xi \in [a,b].
\end{equation}

Now, one can easily write that
\begin{equation}\label{t1}
\| \Phi (x) - \Phi_n(x) \|_{\infty} \leq \|\Phi (x) -S_n(x)
\|_{\infty} + \| S_n(x)- \Phi_n(x)\|_{\infty} \leq \|R_n(x)
\|_{\infty} + \| S_n(x)- \Phi_n(x)\|_{\infty}.
\end{equation}
On the other hand we have
\begin{equation}\label{t2}
\left| S_n(x)- \Phi_n(x) \right| = \left| \sum_{i=0}^n \left(
\widetilde{\Phi}^{(i)}(z) - \Phi^{(i)}(z) \right).  \frac{(x-z)^i}{i
!} \right| \leq |E_i \ell | \leq \| E_i \|_{\infty}  . \| \ell
\|_{\infty} \leq C \| E_i \|_{\infty},
\end{equation}
where $$E_i = \left(e_0(z), e_1(z),...,e_i(z),...,e_n(z)\right),~~~
\ell = \left(\ell_0(z),
\ell_1(z),...,\ell_i(z),...,\ell_n(z)\right),$$ and
$$
e_i(z) = \widetilde{\Phi}^{(i)}(z) - \Phi^{(i)}(z),~~~\ell_i(z) =
\frac{(x-z)^i}{i !}.
$$
By substitute Eqs. (\ref{t0}) and (\ref{t2}) into Eq. (\ref{t1}),
the bound of absolute error function can be obtained as follows
$$
\| \Phi (x) - \Phi_n(x) \|_{\infty} \leq \frac{M}{(i+1)!}
\widetilde{\Phi}^{(i+1)}(\xi) + C \max |e_i(z)| ,~~~\xi \in
[a,b].\blacksquare
$$

\section{Numerical illustrations}
In this section, several examples of generalized Abel's integral
equations of the first and second kinds are solved by using
presented method \cite{d7,a1,a2}. Also, the absolute error for
different values of $x$ and the maximum error are presented in some
tables. The Mathematica 8 was applied to numerical calculations.\\

\noindent \textbf{Example 1:} Consider the following Abel's integral
equation of the first kind \cite{a1,a2}
\begin{equation}\label{4-1}
\frac{2}{3} \pi x^3 = \int_0^x \frac{\Phi(t)}{\sqrt{x^2-t^2}} dt,
\end{equation}
where the exact solution is $\Phi(x) = \pi x^3$. For solving this
integral equation by using the Taylor polynomials the following
algorithm is presented.

\emph{\textbf{\underline{Algorithm}}}

1- Substitute the Taylor polynomial (\ref{2-1}) in Eq. (\ref{4-1})
for arbitrary value $n=5$.
\begin{equation}\label{4-2}
\frac{2}{3} \pi x^3 = \sum_{j=0}^5 \frac{1}{j!} \int_0^x \frac{
\Phi^{(j)}(z) (t-z)^j}{\sqrt{x^2-t^2}} dt.
\end{equation}

2- Substitute the collocation point $x_i$ from Eq. (\ref{2-3}) in
Eq. (\ref{4-2}) as
\begin{equation}\label{4-3}
\frac{2}{3} \pi x_i^3 = \sum_{j=0}^5 \frac{1}{j!} \int_0^{x_i}
\frac{ \Phi^{(j)}(z) (t-z)^j}{\sqrt{x_i^2-t^2}} dt, ~~~~~~~~~~ 0
\leq i \leq n.
\end{equation}

3- Construct the matrix form $AX= G$ where
$$
A= \left[
     \begin{array}{cccccc}
 0& 0& 0& 0& 0& 0\\
 1.5708& 0.2& 0.015708& 0.000888889& 0.0000392699&
 1.42222\times10^{-6}\\
 1.5708& 0.4& 0.0628319& 0.00711111& 0.000628319& 0.0000455111\\
 1.5708& 0.6& 0.141372& 0.024& 0.00318086& 0.0003456\\
 1.5708& 0.8& 0.251327& 0.0568889& 0.0100531& 0.00145636\\
 1.5708& 1& 0.392699& 0.111111& 0.0245437& 0.00444444\\
     \end{array}
   \right],
$$
and
$$
G= \left[
  \begin{array}{c}
    0\\
 0.0167552\\
 0.134041\\
 0.452389\\
 1.07233\\
 2.0944\\
  \end{array}
\right].
$$

4- Determine the unknowns $X$ and calculate the approximate solution
$\Phi_n(x)$
$$
\Phi_5(x) = 0.000211964- 0.00380121 x + 0.0198716 x^2 +
 3.09737 x^3 + 0.0441592 x^4 - 0.0162574 x^5.
$$

The absolute error for $n=5,7,9$ and different values of $x$ is
shown in Table 1 and the maximum errors are demonstrate in Table 2.

\begin{table}[ht]
\caption{Numerical results of Example 1 for $n=5,7,9$.}\label{1}
\centering
\begin{tabular}{||c||c||c||c||c||c||}
  \hline
$x$ & exact     & $e_5 = \mid \Phi(x) -\Phi_{5}(x)   \mid$   &  $e_7 = \mid \Phi(x) -\Phi_{7}(x)   \mid$        &  $e_9 = \mid \Phi(x) -\Phi_{9}(x)   \mid$    \\
    \hline
  0.0    &0            &$0.000211964$                     &$1.17775\times10^{-8}$                          & $1.37683\times10^{-13} $ \\
  0.2  &  0.0251327         &$0.0000417215 $              &$8.60494\times10^{-10}$                         & $ 1.50574\times10^{-15} $      \\
  0.4  & 0.201062       &$ 4.84617\times10^{-6}$              &$1.33074\times10^{-10}$                     & $7.49401\times10^{-16} $ \\
  0.6  &0.678584         &$7.69255\times10^{-6}$                &$7.99231\times10^{-11}$                   &$0$  \\
  0.8  &1.6085        &$8.44368\times10^{-6}$              &$2.03972\times10^{-10}$                        & $4.44089\times10^{-16} $\\
  1.0  & 3.14159           &$0.0000361083$                &$1.60657\times10^{-9}$                          & $9.32587\times10^{-15} $\\
  \hline
  \end{tabular}
 \end{table}

\begin{table}[ht]
\caption{The maximum errors of Example 1.}\label{1} \centering
\begin{tabular}{||c||c||c||}
  \hline
~~~$n=5$~~~&~~~ $n=7$ ~~~   &~~~ $n=9$ ~~~    \\
    \hline
  0.000211964 & $1.17775\times10^{-8}$  & $1.37683\times10^{-13} $  \\
   \hline
  \end{tabular}
\end{table}

\noindent \textbf{Example 2:} Consider the generalized Abel's
integral equation of the first kind \cite{a1}:
\begin{equation}\label{4-4}
\frac{ 4 }{3} \sin^{\frac{3}{4}}(x) = \int_0^x
\frac{\Phi(t)}{(\sin(x) - \sin(t))^{\frac{1}{4}}} dt,
\end{equation}
with the exact solution $ \Phi(x) = \cos x$. By using the presented
method for $n=9$ the approximate solution is obtained as follows
\begin{equation}\label{4-5}
\begin{array}{l}
\Phi_{9}(x) = 0.9999999996930083+ 1.3682136756898444\times 10^{-8}
x -
0.5000001933589282 x^2 \\
~~~~~~~+ 1.3397660985165968\times 10^{-6}  x^3 +
0.041661365627426505 x^4 + \cdots~ .
\end{array}
 \end{equation}
Comparison between exact and approximate solutions for different
values of $x$ and $n=5,7,9$ are presented in Table 3 and the maximum
errors are demonstrate in Table 4.

\begin{table}[ht]
\caption{Numerical results of Example 2 for $n=5,7,9$.}\label{1}
\centering
\begin{tabular}{||c||c||c||c||c||c||}
  \hline
$x$ & exact     & $e_5 = \mid \Phi(x) -\Phi_{5}(x)   \mid$   &  $e_7 = \mid \Phi(x) -\Phi_{7}(x)   \mid$        &  $e_9 = \mid \Phi(x) -\Phi_{9}(x)   \mid$    \\
    \hline
  0.0    &1            &$0.0000590935$                     &$1.82477\times10^{-7}$                          & $3.06992\times10^{-10} $ \\
  0.2  &0.980067         &$0.0000122884 $              &$7.76028\times10^{-9}$                         & $3.44569\times10^{-12} $      \\
  0.4  & 0.921061      &$ 3.76895\times10^{-6}$              &$2.92998\times10^{-9}$                     & $4.82949\times10^{-13} $ \\
  0.6  &0.825336         &$3.85864\times10^{-6}$                &$3.08999\times10^{-9}$                   &$9.76565\times10^{-13} $  \\
  0.8  &0.696707        &$6.42239\times10^{-6}$              &$3.32207\times10^{-9}$                        & $7.80758\times10^{-12} $\\
  1.0  & 0.540302          &$0.000028944$                &$7.77184\times10^{-8}$                          & $1.45912\times10^{-10} $\\
  \hline
  \end{tabular}
 \end{table}

\begin{table}[ht]
\caption{The maximum errors of Example 2.}\label{1} \centering
\begin{tabular}{||c||c||c||}
  \hline
~~~$n=5$~~~&~~~ $n=7$ ~~~   &~~~ $n=9$ ~~~    \\
    \hline
  $0.0000590935$ & $1.82477\times10^{-7}$   & $3.06992\times10^{-10} $  \\
   \hline
  \end{tabular}
\end{table}

\noindent \textbf{Example 3:} In this example, the following first
kind Abel's integral equation \cite{a1} is considered:
\begin{equation}\label{4-6}
 \frac{6}{5} (e^x-1)^{\frac{5}{6}}= \int_0^x \frac{\Phi(t)}{(e^x-e^t)^{\frac{1}{6}}} dt,
 \end{equation}
with exact solution $\Phi(x)= e^x$. By using presented method, we
have following approximate solution
\begin{equation}\label{4-7}
\Phi_9(x) = 1.000000012857142 + 0.9999998135066253 x +
0.5000012303151726 x^2 + \cdots~ .
 \end{equation}

In Table 5 and 6, the absolute errors and maximum values of absolute
error are presented for $n=5,7,9$.

\begin{table}[ht]
\caption{Numerical results of Example 3 for $n=5,7,9$.}\label{1}
\centering
\begin{tabular}{||c||c||c||c||c||c||}
  \hline
$x$ & exact     & $e_5 = \mid \Phi(x) -\Phi_{5}(x)   \mid$   &  $e_7 = \mid \Phi(x) -\Phi_{7}(x)   \mid$        &  $e_9 = \mid \Phi(x) -\Phi_{9}(x)   \mid$    \\
    \hline
  0.0    &1            &$0.00001384$                     &$2.06582\times10^{-7}$                          & $1.28571\times10^{-8} $ \\
  0.2  &1.2214        &$0.0000122884 $              &$4.8359\times10^{-9}$                         & $5.60856\times10^{-10} $      \\
  0.4  &1.49182     &$ 5.25742\times10^{-6}$              &$3.48525\times10^{-9}$                     & $5.87789\times10^{-10} $ \\
  0.6  &1.82212         &$5.33178\times10^{-6}$                &$4.64466\times10^{-9}$                   &$2.37805\times10^{-10} $  \\
  0.8  &2.22554        &$9.82404\times10^{-6}$              &$3.06887\times10^{-9}$                        & $3.93271\times10^{-10} $\\
  1.0  &2.71828          &$0.0000471546$                &$1.27079\times10^{-7}$                          & $1.29719\times10^{-9} $\\
  \hline
  \end{tabular}
 \end{table}

\begin{table}[ht]
\caption{The maximum errors of Example 3.}\label{1} \centering
\begin{tabular}{||c||c||c||}
  \hline
~~~$n=5$~~~&~~~ $n=7$ ~~~   &~~~ $n=9$ ~~~    \\
    \hline
  $0.0000471546$ & $2.06582\times10^{-7}$   &$1.28571\times10^{-8} $  \\
   \hline
  \end{tabular}
\end{table}

\noindent \textbf{Example 4:} In this example, the second kind
Abel's integral equation \cite{a1,d7} is considered as follows
$$
\Phi(x) = x^2 + \frac{16}{15} x^{\frac{5}{2}} - \int_0^x
\frac{\Phi(t)}{(x-t)^{\frac{1}{2}}} dt,
$$
with the exact solution $\Phi(x) = x^2$. Presented method for
solving this example is exact for $n=5$.

\noindent \textbf{Example 5:} Let us consider the following integral
equation \cite{a1}:
$$
\Phi(x) = 1-2x-\frac{32}{21} x^{\frac{7}{4}} + \frac{4}{3}
x^{\frac{3}{4}} - \int_0^x \frac{\Phi(t)}{(x-t)^{\frac{1}{4}}} dt,
$$
with the exact solution $\Phi(x) = 1-2x$. In this example, the
Taylor-collocation method for $n=3$ is exact.

\section*{Acknowledgements} This work has been supported by Islamic Azad
University, Hamedan branch, research grand.

\section{Conclusions}
In this paper, Taylor polynomials to estimate the generalized form
of Abel's integral equations of the first and second kinds were
used. Also, bound of absolute error function was illustrated.
Several applicable examples were solved by using Taylor-collocation
method. Presented tables and obtained numerical results show the
efficiency and accuracy of method.



\begin{thebibliography}{99}

\bibitem{5} CT.H. Baker, The numerical treatment of integral equations.
Oxford: Clarendon Press; 1977.


\bibitem{b7} H. Brunner, 1896-1996: One hundred years of Volterra integral
equation of the first kind, Appl. Numer. Math. 24 (1997) 83-93.

\bibitem{b4} C. J. Cremers, R.C. Birkebak, Application of the Abel
integral equation to spectroscopic data, Appl. Opt. 5 (1966)
1057-1064.

\bibitem{b3} S. De, B. N. Mandal, A. Chakrabarti, Use of Abel's integral
equations in water wave scattering by two surface-piercing barriers,
Wave Motion 47 (2010) 279-288.

\bibitem{8} L.M. Delves, JL. Mohamed, Computational methods for integral
equations. UK: Cambridge University Press; 1985.

\bibitem{d2} M.A. Fariborzi Araghi, H. Daei Kasmaei, Numerical Solution of the Second Kind Singular Volterra Integral
Equations By Modified Navot-Simpson's Quadrature, Int. J. Open
Problems Compt. Math. 1 (2008).

\bibitem{col2}M.A. Fariborzi Araghi, Gh. Kazemi
Gelian, Numerical solution of nonlinear Hammerstein integral
equations via Sinc-collocation method based on double exponential
transformation, Mathematical Sciences, 7 (2013) 1-7.

\bibitem{col1} M.A. Fariborzi Araghi, S. Noeiaghdam,
Fibonacci-regularization method for solving Cauchy integral
equations of the first kind, Ain Shams Eng J. 8 (2017) 363-369.

\bibitem{myy} M.A. Fariborzi Araghi, S. Noeiaghdam, A novel technique based on the
homotopy analysis method to solve the first kind Cauchy integral
equations arising in the theory of airfoils, {Journal of
Interpolation and Approximation in Scientific Computing } 2016 (1)
 (2016) 1-13.

\bibitem{me1} M.A. Fariborzi Araghi, S. Noeiaghdam, Homotopy analysis transform
method for solving generalized Abel's fuzzy integral equations of
the first kind, IEEE (2016). DOI: 10.1109/CFIS.2015.7391645

\bibitem{me2} M.A. Fariborzi Araghi, S. Noeiaghdam, Homotopy regularization method
to solve the singular Volterra integral equations of the first kind,
Jordan Journal of Mathematics and Statistics (JJMS) 11(1) 2018,
1-12.


\bibitem{d1}M.A. Fariborzi Araghi, S. Yazdani, A Method to Approximate Solution
of the First kind Abel Integral Equation using Navot's Quadrature
and Simpson's rule, International Journal of Industrial Mathematics,
 1 (2009) 1-11.


\bibitem{b5} R. Gorenflo, S. Vessella, Abel Integral Equations: Analysis
and Applications, Lecture Notes Math., 1461, Springer, Berlin, 1991.


\bibitem{e2} P. Huabsomboon, B. Novaprateep, H. Kaneko, On Taylor-series expansion methods for the second kind integral
equations, J. Comput. Appl. Math. 234 (2010) 1466-1472.


\bibitem{d6} L. Huang, Y. Huang, X. F. Li, Approximate solution of Abel integral equation, Comput. Math. Appl. 56 (2008) 1748-1757.


\bibitem{d4} S. Jahanshahi, E. Babolian, Delfim F.M. Torres, A. Vahidi ,Solving Abel integral equations of first kind via
fractional calculus, Journal of King Saud University-Science (27)
2015 161-167.

\bibitem{e5} R.P. Kanwal, K.C. Liu, A Taylor expansion approach for solving
integral equations, Int. J. Math. Educ. Sci. Technol. 20 (1989)
411-414.

\bibitem{b2} O. Knill, R. Dgani, M. Vogel, A new approach to Abel's
integral operator and its application to stellar winds, Astronom.
Astrophys. 274 (1993) 1002-1008.


\bibitem{b6} L. E. Kosarev, Applications of integral equations of the first
kind in experiment physics., Comput. Phys. Commun. 20 (1980) 69-75.


\bibitem{b11} P.K. Lamm, L. Elden, Numerical solution of first-kind Volterra
equations by sequential Tikhonov regularization, SIAM J. Numer.
Anal. 34 (4) (1997) 1432-1450.

\bibitem{col4}K. Maleknejad, K. Nedaiasl, Application of Sinc-collocation method
for solving a class of nonlinear Fredholm integral equations,
Comput. Math. Appl. 62 (2011) 3292-3303.

\bibitem{e3} K. Maleknejad, N. Aghazadeh, M. Rabbani, Numerical omputational
solution of second kind Fredholm integral equations system by using
a Taylor-series expansion method, Appl. Math. Comput. 175 (2006)
1229-1234.

\bibitem{e9} K. Maleknejad, Y. Mahmoudi, Taylor polynomial solution of
high-order nonlinear Volterra-Fredholm integro-differential
equations, Appl. Math. Comput. 145 (2003) 641-653.



\bibitem{d7} S. Noeiaghdam, E. Zarei, H. Barzegar Kelishami, Homotopy
analysis transform method for solving Abel's integral equations of
the first kind, {Ain Shams Eng. J.}  7 (2016) 483-495.


\bibitem{b9} A. D. Polyanin, A.V. Manzhirov, Handbook of Integral
Equations, CRC Press, Boca Raton, 1998.

\bibitem{e6} M. Sezer, Taylor polynomial solution of Volterra integral
equations, Int. J. Math. Educ. Sci. Technol. 25 (1994) 625-633.


\bibitem{d3}S. Sohrabi, Comparison Chebyshev wavelets method with BPFs method for solving
Abel's integral equation , Ain Shams Engineering Journal (2) 2011
249-254.


\bibitem{b10} H.J. Teriele, Collocation method for weakly singular second kind
Volterra integral equations with non-smooth solution, IMA J. Numer.
Anal. 2 (4) (1982) 437-449.


\bibitem{b1} J. Walton, Systems of generalized Abel integral equations with
applications to simultaneous dual relations, SIAM J. Math. Anal. 10
 (1979) 808-822.

\bibitem{e1} K. Wang, Q. Wang, Taylor polynomial method and error estimation for a kind of mixed
Volterra-Fredholm integral equations, Appl. Math. Comput. 229 (2014)
53-59.

\bibitem{e4} K. Wang, Q. Wang, Taylor collocation method and convergence analysis for the
Volterra-Fredholm integral equations, J. Comput. Appl. Math. 260
(2014) 294-300.

\bibitem{a1} A.M. Wazwaz, Linear and nonlinear integral equations: methods
and applications. Berlin: Higher Education, Beijing, and Springer;
2011.

\bibitem{a2} A.M. Wazwaz, A First Course in Integral Equations, World
Scientific Publishing, Singapore, 1997.


\bibitem{a4} A.M. Wazwaz, M.S. Mehanna, The combined Laplace-Adomian method
for handling singular integral equation of heat transfer, Int. J.
Nonlinear Sci. 10 (2010) 248-252.

\bibitem{e7} S. Yalcinbas, Taylor polynomial solutions of nonlinear
Volterra-Fredholm integral equations, Appl. Math. Comput. 127 (2002)
195-206.


\bibitem{e8} S. Yalcinbas, M. Sezer, The approximate solution of high-order
linear Volterra-Fredholm integro-differential equations in terms of
Taylor polynomials, Appl. Math. Comput. 112 (2000) 291-303.

\bibitem{d5} C. Yang, An efficient numerical method for solving Abel integral equation, Appl. Math. Comput. 227 (2014) 656-661.

\bibitem{a5} N. Zeilon, Sur quelques points de la theorie de
l'equationintegraled'Abel, Arkiv. Mat. Astr. Fysik. 18 (1924) 1-19.


\end{thebibliography}
\end{document}